\documentclass[leqno]{article}
\usepackage{amsmath}
\usepackage{amssymb}
\def\Re{\mathop{\mathrm{Re}}}
\def\ro{{\varrho}}
\def\eps{\varepsilon}

\begin{document}

\title{\sc On certain arithmetic functions involving exponential divisors, II.}
\author{{\bf L\'aszl\'o T\'oth} (P\'ecs, Hungary)}
\date{}
\maketitle

\centerline{\it Dedicated to the memory of Professor M. V. Subbarao}
\vskip2mm

\centerline{Annales Univ. Sci. Budapest., Sect. Comp. {\bf 27} (2007), 155-166}

\vskip2mm {\bf 1. Introduction}

\vskip1mm Let $n>1$ be an integer of canonical form $n=p_1^{a_1}\cdots p_r^{a_r}$.
The integer $d$ is called an exponential divisor (e-divisor) of $n$ if
$d=p_1^{b_1}\cdots p_r^{b_r}$, where $b_1 \mid a_1$, ..., $b_r\mid a_r$, notation:
$d\mid_e n$. By convention $1\mid_e 1$. The integer $n>1$ is called exponentially
squarefree (e-squarefree) if all the exponents $a_1,...,a_r$ are squarefree. The
integer $1$ is also considered to be e-squarefree.

The exponential convolution (e-convolution) of arithmetic functions
is defined by
\[
(f\odot g)(n)=\sum_{b_1c_1=a_1} \dots \sum_{b_rc_r=a_r} f(p_1^{b_1}\cdots
p_r^{b_r}) g(p_1^{c_1}\cdots p_r^{c_r}),
\]
where $n=p_1^{a_1}\cdots p_r^{a_r}$.

These notions were introduced by M. V. Subbarao \cite{Su72}. The e-convolution
$\odot$ is commutative, associative and has the identity element $\mu^2$, where
$\mu$ is the M\"obius function. Furthermore, a function $f$ has an inverse with
respect to $\odot$ iff $f(1)\ne 0$ and $f(p_1\cdots p_s)\ne 0$ for any distinct
primes $p_1,...,p_s$.

The inverse with respect to $\odot$ of the constant $1$ function is called the
exponential analogue of the M\"obius function and it is denoted by $\mu^{(e)}$.
Hence for every $n\ge 1$,
\[
\sum_{d\mid_e n} \mu^{(e)}(d)=\mu^2(n).
\]

Here $\mu^{(e)}(1)=1$ and for  $n=p_1^{a_1}\cdots p_r^{a_r}>1$,
\[
\mu^{(e)}(n)=\mu(a_1)\cdots \mu(a_r).
\]

Observe that $|\mu^{(e)}(n)|=1$ or $0$, according as $n$ is e-squarefree or not.
For properties and generalizations of the e-convolution see \cite{Su72},
\cite{HR97}.

Other arithmetic functions regarding e-divisors, for example the number and the sum
of e-divisors of $n$ were investigated by several authors, see the references given
in the first part \cite{To2004} of the present paper, devoted to the study of
functions involving the greatest common exponential divisor of integers.

An asymptotic formula for $\sum_{n\le x} |\mu^{(e)}(n)|$ was established by M. V.
Subbarao \cite{Su72}, improved by J. Wu \cite{Wu95}, see also Part I. of the
present paper. We show that the corresponding error term can further be improved on
the assumption of the Riemann hypothesis (RH), see Theorem 3.

In Theorem 2 we give a formula for $\sum_{n\le x} \mu^{(e)}(n)$
without and with assuming RH. As far as we know there is no such
result in the literature. We show that the error terms depend on
estimates for the number of squarefree integers $\le x$.

Consider now the exponential squarefree exponential divisors (e-squarefree
e-divisors) of $n$. Here $d=p_1^{b_1}\cdots p_r^{a_r}$ is an e-squarefree e-divisor
of $n=p_1^{a_1}\cdots p_r^{a_r}>1$, if $b_1 \mid a_1,..., b_r\mid a_r$ and
$b_1,...,b_r$ are squarefree. Note that the integer $1$ is e-squarefree and it is
not an e-divisor of $n>1$.

We introduce the functions $t^{(e)}$ and $\kappa^{(e)}$, where $t^{(e)}(n)$ and
$\kappa^{(e)}(n)$ denote the number of e-squarefree e-divisors of $n$ and the
maximal e-squarefree e-divisor of $n$, respectively. These are the exponential
analogues of the functions representing the number of squarefree divisors of $n$
(i.e. $\theta(n)=2^{\omega(n)}$, where $\omega(n)=r$) and the maximal squarefree
divisor of $n$ (the squarefree kernel $\kappa(n)=\prod_{p\mid n} p$), respectively.

The functions $t^{(e)}$ and $\kappa^{(e)}$ are multiplicative and for
$n=p_1^{a_1}\cdots p_r^{a_r}>1$,
\[
t^{(e)}(n)=2^{\omega(a_1)}\cdots 2^{\omega(a_r)},
\]
\[
\kappa^{(e)}(n)=p_1^{\kappa(a_1)}\cdots p_r^{\kappa(a_r)}.
\]

Asymptotic properties of the functions $t^{(e)}(n)$ and $\kappa^{(e)}(n)$ are given
in Theorems 4, 5 and 7.

\vskip2mm {\bf 2. Results}

\vskip1mm The function $\mu^{(e)}$ is multiplicative and $\mu^{(e)}(p^a)=\mu(a)$
for every prime power $p^a$. Hence $\mu^{(e)}(n) \in \{-1,0,1\}$ for every $n\ge 1$
and for every prime $p$, $\mu^{(e)}(p)=1$, $\mu^{(e)}(p^2)=-1$,
$\mu^{(e)}(p^3)=-1$, $\mu^{(e)}(p^4)=0$,... .

According to a well-known result of H. Delange, cf. \cite{El79}, Ch.
6, the function $\mu^{(e)}$ has a non-zero mean value given by
\[
m(\mu^{(e)}) =\prod_p \left(1 + \sum_{a=2}^{\infty} \frac{\mu(a)-\mu(a-1)}{p^a}
\right).
\]

An asymptotic formula for $\mu^{(e)}$ can be obtained from the
following general result, which may be known.

\vskip1mm {\bf Theorem 1.} {\it Let $f$ be a complex valued
multiplicative function such that $|f(n)|\le 1$ for every $n\ge 1$
and $f(p)=1$ for every prime $p$. Then
\[
\sum_{n\le x} f(n)= m(f) x +O(x^{1/2} \log x),
\]
where
\[
m(f)=\prod_p \left(1 + \sum_{a=2}^{\infty}
\frac{f(p^a)-f(p^{a-1})}{p^a} \right).
\]
is the mean value of $f$.}

\vskip1mm Theorem 1 applies also for the multiplicative functions $f=\mu^{*(e)}$
and $f=F$, where $\mu^{*(e)}(p^a)=\mu^*(a)=(-1)^{\omega(a)}$ representing the
unitary exponential M\"obius function, cf. \cite{HR97}, and
$F(p^a)=\lambda(a)=(-1)^{\Omega(a)}$ the Liouville function, with $\Omega(a)$
denoting the number of prime power divisors of $a$.

We prove for $\mu^{(e)}$ the following more precise result.

\vskip1mm {\bf Theorem 2.} {\it (i) The Dirichlet series of
$\mu^{(e)}$ is of form
\[
\sum_{n=1}^{\infty} \frac{\mu^{(e)}(n)}{n^s} = \frac{\zeta(s)}{\zeta^2(2s)} U(s),
\quad \Re s > 1,
\]
where $U(s):=\sum_{n=1}^{\infty} \frac{u(n)}{n^s}$ is absolutely convergent for
$\Re s > 1/5$.

(ii)
\[
\sum_{n\le x} \mu^{(e)}(n)= m(\mu^{(e)})x+ O(x^{1/2}\exp(-c (\log x)^{\Delta}),
\]
where $\Delta<9/25= 0,36$ and $c>0$ are constants.

(iii) Assume RH. Let $1/4<r<1/3$ be an exponent such that $D(x):=
\sum_{n\le x} \mu^2(n)-x/\zeta(2)= O(x^{r+\eps})$ for every $\eps
>0$. Then the error term in (ii) is $O(x^{(2-r)/(5-4r)+\eps})$ for
every $\eps >0$. }

\vskip1mm
The best known value -- to our knowledge -- of $r$ is $r=17/54\approx
0,314814$, obtained in \cite{Ji93}, therefore the error term in (ii), assuming RH,
is $O(x^{91/202+\eps})$ for every $\eps >0$, where $91/202\approx 0,450495$.

\vskip1mm {\bf Theorem 3.} {\it If RH is true, then
\[
\sum_{n\le x} |\mu^{(e)}(n)| = \prod_p \left(1 + \sum_{a=4}^{\infty}
\frac{\mu^2(a)-\mu^2(a-1)}{p^a} \right)x+ O(x^{1/5+\varepsilon}),
\]
for every $\varepsilon >0$.}

\vskip1mm The function $t^{(e)}$ is multiplicative and $t^{(e)}(p^a)=2^{\omega(a)}$
for every prime power $p^a$. Here for every prime $p$, $t^{(e)}(p)=1$,
$t^{(e)}(p^2)=t^{(e)}(p^3)=t^{(e)}(p^4)=t^{(e)}(p^5)=2, t^{(e)}(p^6)=4,...$ .

\vskip1mm {\bf Theorem 4.} {\it (i) The Dirichlet series of $t^{(e)}$ is of form
$$
\sum_{n=1}^{\infty} \frac{t^{(e)}(n)}{n^s} =\zeta(s)\zeta(2s) V(s), \quad \Re s>1,
$$
where $V(s)=\sum_{n=1}^{\infty} \frac{v(n)}{n^s}$ is absolutely convergent for $\Re
s >1/4$.

(ii)
\[
\sum_{n\le x} t^{(e)}(n)= C_1x +C_2 x^{1/2} + O(x^{1/4+\eps}),
\]
for every $\eps >0$, where $C_1, C_2$ are constants given by
\[
C_1:=\prod_p \left(1 + \frac1{p^2} + \sum_{a=6}^{\infty}
\frac{2^{\omega(a)}-2^{\omega(a-1)}}{p^a}\right),
\]
\[
C_2:=\zeta(1/2) \prod_p \left(1 + \sum_{a=4}^{\infty}
\frac{2^{\omega(a)}-2^{\omega(a-1)}-2^{\omega(a-2)}+
2^{\omega(a-3)}}{p^{a/2}}\right).
\]}

\vskip1mm {\bf Theorem 5.} {\it
\[
\limsup_{n\to \infty} \frac{\log t^{(e)}(n) \log \log n}{\log n} = \frac1{2} \log
2.
\]}

\vskip1mm The function $\kappa^{(e)}$ is multiplicative and
$\kappa^{(e)}(p^a)=p^{\kappa(a)}$ for every prime power $p^a$. Hence for every
prime $p$, $\kappa^{(e)}(p)=p$, $\kappa^{(e)}(p^2)=p^2$, $\kappa^{(e)}(p^3)=p^3$,
$\kappa^{(e)}(p^4)=p^2$,... .

To obtain an asymptotic formula for $\kappa^{(e)}$ we use the
following general theorem, of which parts (i) and (ii) are a variant
of a result given in \cite{SiSu82} and cited in the first part
\cite{To2004} of this paper.

\vskip1mm {\bf Theorem 6.} {\it Let $k\ge 2$ be a fixed integer and
$f$ be a complex valued multiplicative arithmetic function
satisfying

(a) $f(p)=f(p^2)=...=f(p^{k-1})=1$ for every prime $p$,

(b) there exists $K>0$ such that $|f(p^a)|\le K$ for every prime
power $p^a$ with $a\ge k+1$,

(c) there exist $M>0$ and $\beta \ge 1/(k+1)$ such that $|f(p^k)|\le
Mp^{-\beta}$ for every prime $p$ .

Then

(i) \[ \sum_{n=1}^{\infty} \frac{f(n)}{n^s} =
\frac{\zeta(s)}{\zeta(ks)} W(s), \quad \Re s > 1,
\]
where the Dirichlet series $W(s):=\sum_{n=1}^{\infty} \frac{w(n)}{n^s}$ is
absolutely convergent for $\Re s > 1/(k+1)$.

(ii) \[ \sum_{n\le x} f(n)=C_f x +O(x^{1/k}\delta(x)),
\]
where
\[
C_f:=\prod_p \left(1+ \sum_{a=k}^{\infty}
\frac{f(p^a)-f(p^{a-1})}{p^a}\right)
\]
and
\[
\delta(x)=\delta_A(x):=\exp(-A(\log x)^{3/5} (\log \log x)^{-1/5}),
\]
$A$ being a positive constant.

(iii) If RH is true, then the error term is
$O(x^{1/(k+1)+\varepsilon})$ for every $\varepsilon >0$.}

\vskip1mm {\bf Theorem 7.} {\it
\[
\sum_{n\le x} \kappa^{(e)}(n)=\frac1{2} \prod_p \left(1 + \sum_{a=4}^{\infty}
\frac{p^{\kappa(a)}-p^{1+\kappa(a-1)}}{p^a} \right) x^2 + O(x^{5/4}\delta(x)).
\]

If RH is true, then the error term is $O(x^{6/5+\varepsilon})$ for every
$\varepsilon >0$.}

\vskip2mm {\bf 3. Proofs}

\vskip1mm {\bf Proof of Theorem 1.} Let $g=f*\mu$ in terms of the
Dirichlet convolution. Then $g$ is multiplicative, $g(p)=f(p)-1=0$,
$g(p^a)=f(p^a)-f(p^{a-1})$ and $|g(p^a)|\le |f(p^a)|+|f(p^{a-1})|\le
2$ for every prime $p$ and every $a\ge 2$. Therefore $|g(n)|\le
\ell(n) 2^{\omega(n)}$ for every $n\ge 1$, where $\ell(n)$ is the
characteristic function of the squarefull integers and we have
\[
\sum_{n\le x} f(n)=\sum_{de\le x} g(d)= \sum_{d\le x} g(d)
\left(\frac{x}{d}+O(1)\right)= x\sum_{d\le x}
\frac{g(d)}{d}+O\left(\sum_{d\le x} |g(d)|\right)=
\]
\[ =x \sum_{d=1}^{\infty} \frac{g(d)}{d} + O\left(x \sum_{d>x}
\frac{\ell(d)2^{\omega(d)}}{d}\right) + O\left(\sum_{d\le x} \ell(d)
2^{\omega(d)}\right).
\]

Here
\[
\ell(n)2^{\omega(n)}=\sum_{d^2e=n} \tau(d)h(e),
\]
where $\tau$ is the divisor function and $h$ is given by
\[
\sum_{n=1}^{\infty} \frac{h(n)}{n^s}=\prod_p
\left(1+\frac{2}{p^{3s}}-\frac{1}{p^{4s}}-\frac{2}{p^{5s}}\right),
\]
absolutely convergent for $\Re s>1/3$, cf. \cite{SS79}. We obtain
\[
\sum_{n\le x} \ell(n)2^{\omega(n)}=\sum_{e\le x} h(e) \sum_{d\le
(x/e)^{1/2}} \tau(d) = \sum_{e\le x} h(e) \ O\left( (x/e)^{1/2}\log
(x/e) \right)= \]
\[
=O\left(x^{1/2}\log x \sum_{e\le x} |h(e)|e^{-1/2}\right)=
O\left(x^{1/2}\log x \right),
\]
and by partial summation,
\[
\sum_{n>x} \frac{\ell(n)2^{\omega(n)}}{n}= O\left(x^{-1/2}\log x
\right),
\]
which finishes the proof.

\vskip1mm {\bf Proof of Theorem 2.} (i) Let $\mu_2(n)=\mu(m)$ or $0$, according as
$n=m^2$ or not, and let $E_2(n)=1$ or $0$, according as $n=m^2$ or not. The given
equality is verified for $\mu^{(e)}=\mu^2*\mu_2*u$, equivalent to
$u=\mu^{(e)}*\lambda* E_2$, in terms of the Dirichlet convolution, where $\lambda$
is the Liouville function. It is easy to check that $u(p)=u(p^2)=u(p^3)=u(p^4)=0$,
$|(\lambda* E_2)(p^a)|\le a$ for every prime power $p^a$ with $a\ge 1$, hence
$|u(p^b)|\le 1+ \sum_{a=1}^b |(\lambda* E_2)(p^a)|< b^2$ for every prime power
$p^b$ with $b\ge 5$. We obtain that the Dirichlet series of the function $u$ is
absolutely convergent for $\Re s> 1/5$.

(ii) According to (i), $\sum_{n\le x} \mu^{(e)}(n) =\sum_{n\le x}
u(n)S(x/n)$, where
\[ S(x):=\sum_{nd^2\le x} \mu^2(n)\mu(d).\]

We first estimate the sum $S(x)$. Let $\varrho=\varrho(x)$ such that $0<\varrho<1$
to be defined later. If $nd^2\le x$, then both $n>\varrho^{-2}$ and $d>\ro
\sqrt{x}$ can not hold good in the same time, therefore
\[
S(x)=\sum_{\substack{nd^2\le x \\ d\le \ro \sqrt{x}}} \mu^2(n)\mu(d)+
\sum_{\substack{nd^2\le x \\ n\le \ro^{-2}}} \mu^2(n) \mu(d)- \sum_{\substack{d\le
\ro \sqrt{x} \\ n\le \ro^{-2}}} \mu^2(n)\mu(d) = S_1(x) + S_2(x) - S_3(x),
\]
say. We use the following estimates of A. Walfisz \cite{Wa63}:
\[
M(x):= \sum_{n\le x} \mu(n) =O(x\delta(x)), \quad E(x):= \sum_{n\le
x} \mu^2(n) =\frac{x}{\zeta(2)}+ O(x^{1/2} \delta(x)).
\]

Note that $\delta(x)$, defined in Section 2, is decreasing and $x^{\eps}\delta(x)$
is increasing for every $\eps>0$. By partial summation,
\[
R(x):=\sum_{n>x} \frac{\mu(n)}{n^2} =O(x^{-1}\delta(x)).
\]

Here
\[
S_1(x)=\sum_{d\le \ro \sqrt{x}} \mu(d) E(x/d^2)= \frac{x}{\zeta(2)}\sum_{d\le \ro
\sqrt{x}} \frac{\mu(d)}{d^2}+ O\left(x^{1/2}\sum_{d\le \ro \sqrt{x}}
\frac{\delta(x/d^2)}{d}\right)= \]
\[= \frac{x}{\zeta(2)} \left(\frac1{\zeta(2)}- R(\ro
\sqrt{x})\right) + O\left(x^{1/2}\delta(\ro^{-2}) \sum_{d\le \sqrt{x}}
\frac1{d}\right)= \] \[= \frac{x}{\zeta^2(2)}+ O\left(\ro^{-1} x^{1/2} \delta(\ro
\sqrt{x}) \right) + O\left(x^{1/2}\delta(\ro^{-2})\log x \right),
\]
\[
S_2(x)=\sum_{n\le \ro^{-2}} \mu^2(n) M((x/n)^{1/2})=O\left(\sum_{n\le \ro^{-2}}
(x/n)^{1/2} \delta((x/n)^{1/2})\right)= O\left(\delta(\ro \sqrt{x}) x^{1/2}
\sum_{n\le \ro^{-2}} \frac1{\sqrt{n}} \right)= \]
\[ = O\left(\ro^{-1} x^{1/2} \delta(\ro \sqrt{x}) \right),
\]
\[
S_3(x)=M(\ro \sqrt{x}) E(\ro^{-2})=O\left(\ro^{-1} x^{1/2} \delta(\ro \sqrt{x})
\right).
\]

We obtain that
\[ S(x)=\frac{x}{\zeta^2(2)}+ O\left(\ro^{-1} x^{1/2} \delta(\ro \sqrt{x})
\right) + O\left(x^{1/2}\delta(1/\ro^2)\log x \right).
\]

Take $\ro = \exp(-(\log x)^{\beta})$, where $0<\beta<1$. Then $\ro
\sqrt{x}=\exp(\frac1{2}(\log x)-(\log x)^{\beta})\ge \exp(\frac1{4}(\log
x))=x^{1/4}$ for sufficiently large $x$. Hence $\delta(\ro \sqrt{x})\le
\delta(x^{1/4})\ll \delta_B(x)$ with a suitable constant $B>0$. For $\beta<3/5$ we
obtain $\ro^{-1} \delta(\ro \sqrt{x})\ll \exp((\log x)^{\beta} - B(\log x)^{3/5}
(\log \log x)^{-1/5})\ll \delta_C(x)$ with a suitable constant $C>0$.

If $\eta<3/5$, then $\delta_A(x)\ll \exp(-A(\log x)^{\eta})$ and obtain that
$\delta(\ro^{-2})\ll \exp(-A(2(\log x)^{\beta})^{\eta}) = \exp(-D(\log x)^{\beta
\eta})$ with a suitable $D>0$, where $\beta \eta< 9/25$.

Therefore,
\[
S(x)= \frac{x}{\zeta^2(2)}+ O\left(x^{1/2}\exp(-c (\log
x)^{\Delta}\right),
\]
where $\Delta<9/25$ and $c>0$ are constants. Now,
\[
\sum_{n\le x} \mu^{(e)}(n)= \sum_{n\le x} u(n) S(x/n)= \sum_{n\le x}
u(n) \left(\frac{x}{\zeta^2(2)n}+ O\left((x/n)^{1/2}\exp(-c (\log
(x/n))^{\Delta}) \right)\right)= \] \[= \frac{x}{\zeta^2(2)}
\sum_{n\le x} \frac{u(n)}{n} +O\left(x^{1/2} \sum_{n\le x}
\frac{|u(n)|}{n^{1/2}} \exp(-c (\log (x/n))^{\Delta})\right),
\]
where, using that $x^{\eps} \exp(-c(\log x)^{\Delta})$  is
increasing for any $\eps>0$, the $O$-term is
\[ O\left(x^{1/2} \sum_{n\le x} \frac{|u(n)|}{n^{1/2}}
(\frac{x}{n})^{-\eps} (\frac{x}{n})^{\eps} \exp(-c (\log(x/n)^{\Delta})) \right)=
O\left( x^{1/2} x^{\eps} \exp(-c(\log x)^{\Delta}) x^{-\eps} \sum_{d\le x}
\frac{|u(n)|}{n^{1/2-\eps}} \right)= \] \[ =O\left(x^{1/2}\exp(-c(\log
x)^{\Delta})\right),
\]
for $1/2- \eps>1/5$. Furthermore,
\[
\sum_{n\le x} \frac{u(n)}{n}=U(1)+ O\left(\sum_{n>x}
\frac{|u(n)|}{n}\right),
\] with $U(1)=\zeta^{-2}(2) m(\mu^{(e)})$ and $\sum_{n>x}
\frac{|u(n)|}{n}= O\left(x^{-3/5}\sum_{n>x}
\frac{|u(n)|}{n^{2/5}}\right) =O(x^{-3/5})$, which finishes the proof
of (ii).

(iii) Assume RH. We use that, see \cite{Ti51},
\[
M(x):= \sum_{n\le x} \mu(n)= O\left(x^{1/2}\omega(x)\right),
\]
where $\omega(x):=\exp(A(\log x)(\log \log x)^{-1})$, $A$ being a
positive constant, which gives by partial summation,
\[
R(x):=\sum_{n>x} \frac{\mu(n)}{n^2} =O(x^{-3/2}\omega(x)).
\]

Suppose that $D(x):=\sum_{n\le x} \mu^2(n)-x/\zeta(2)=O(x^{r+\eps})$
for every $\eps>0$, where $1/4<r<1/3$. Then we obtain by similar
computations that
\[
S_1(x)= \frac{x}{\zeta^2(2)}+ O\left(\ro^{-3/2} x^{1/4} \omega(\ro
\sqrt{x}) \right) + O\left(x^{1/2} \ro^{1-2(r+\eps)} \right),
\]
\[
S_2(x)=O\left(\ro^{-3/2} x^{1/4} \omega(\ro \sqrt{x}) \right),
\qquad S_3(x)= O\left(\ro^{-3/2} x^{1/4} \omega(\ro \sqrt{x})
\right),
\]

Therefore
\[
S(x)= \frac{x}{\zeta^2(2)}+ O\left(\ro^{-3/2} x^{1/4} \omega(\ro
\sqrt{x}) \right) + O\left(x^{1/2} \ro^{1-2(r+\eps)} \right).
\]

Choose $\varrho=x^{-t}$, $t>0$. Then $\ro^{-3/2}
x^{1/4}=x^{(6t+1)/4}$, $\ro \sqrt{x}=x^{1/2-t}<x$, hence $\omega(\ro
\sqrt{x})<\omega(x)\ll x^{\eps}$ for every $\eps >0$ and obtain
\[
S(x)= \frac{x}{\zeta^2(2)}+ O\left(x^{(6t+1)/4+\eps} \right) +
O\left(x^{1/2-t(1-2r)+\eps}\right).
\]

Take $(6t+1)/4 = 1/2-t(1-2r)$, this gives $t=1/(10-8r)$ leading to the common value
$(2-r)/(5-4r)+\eps$ of the exponents.

\vskip1mm {\bf Proof of Theorem 3.} Apply Theorem 6 for
$f(n)=|\mu^{(e)}(n)|$, $k=4$ on the assumption of RH.

\vskip1mm {\bf Proof of Theorem 4.} The proof is similar to the
proof of Theorem 1 of \cite{To2004}, see also \cite{T} for a more
general result of this type.

(i) To obtain the given equality let $f=\mu_2*\mu$, where $\mu_2$ is defined in the
Proof of Theorem 2, and let $v=t^{(e)}*f$. Here both $f$ and $v$ are multiplicative
and it is easy to check that $f(p)=f(p^2)=-1, f(p^3)=1, f(p^a)=0$ for each $a\ge
4$, and $v(p)=v(p^2)=v(p^3)=0$, $v(p^a)=2^{\omega(a)}-2^{\omega(a-1)}-
2^{\omega(a-2)}+ 2^{\omega(a-3)}$ for $a\ge 4$.

(ii) According to (i), $t^{(e)}=v*\tau(1,2,\cdot)$, where
$\tau(1,2,n)=\sum_{ab^2=n} 1$ for which
\[
\sum_{n\le x} \tau(1,2,n) =\zeta(2)x+\zeta(1/2)x^{1/2}+O(x^{1/4}),
\]
cf. \cite{Kr88}, p. 196-199. Therefore,
\[
\sum_{n\le x} t^{(e)}(n) =\sum_{d\le x} v(d) \sum_{e \le x/d} \tau(1,2,e)
\]
and we obtain the above result by usual estimates.

\vskip1mm {\bf Proof of Theorem 5.} We use the following general
result given in \cite{SuSi75}: Let $F$ be a multiplicative function
with $F(p^a)=f(a)$ for every prime power $p^a$, where $f$ is
positive and satisfying $f(n)=O(n^\beta)$ for some fixed $\beta
>0$. Then
\[
\limsup_{n\to \infty} \frac{\log F(n) \log \log n}{\log n} = \sup_{m} \frac{\log
f(m)}{m}.
\]

Take $F(n)=t^{(e)}(n)$, $f(a)=2^{\omega(a)}$. Here $\omega(a)\le
a/2$ and $\frac{\log f(2)}{2}=\frac1{2}\log 2$, which proves the
result.

\vskip1mm {\bf Proof of Theorem 6.} (i), (ii) Take $f=q_k*w$, in terms of the
Dirichlet convolution, where $q_k$ denotes the characteristic function of the
$k$-free integers and use the estimate of A. Walfisz \cite{Wa63},
$$
\sum_{n\le x} q_k(n)=\frac{x}{\zeta(k)}+O(x^{1/k}\delta(x)).
$$
For details cf. \cite{SiSu82}, \cite{T}.

(iii) If RH is true, then the error term of above is $O(x^{1/(k+1)+\varepsilon})$,
according to the result of H. L. Montgomery and R. C. Vaughan \cite{MV81}, and take
into account that $W(s)$ is absolutely convergent for $\Re s> 1/(k+1)$.

\vskip1mm {\bf Proof of Theorem 7.} Apply Theorem 6 for
$f(n)=\kappa^{(e)}(n)/n$, $k=4$, $\beta=2$, where $f(p^4)=1/p^2$.
Then by partial summation we obtain the result.

\vskip2mm {\bf Acknowledgement.} The author is grateful to Professor
Imre K\'atai for valuable suggestions.

\vskip4mm

\noindent{{\bf L\'aszl\'o T\'oth}\\
University of P\'ecs\\
Institute of Mathematics and Informatics\\
Ifj\'us\'ag u. 6\\
7624 P\'ecs, Hungary\\
ltoth@ttk.pte.hu}

\end{document}